\documentclass[12pt]{amsart}
\usepackage{amssymb}

\newtheorem{thm}{Theorem}
\newtheorem{prop}[thm]{Proposition}
\newtheorem{cor}[thm]{Corollary}
\newtheorem{conj}[thm]{Conjecture}

\begin{document}
\title[The special Schubert calculus is real]{The special Schubert calculus
is real}   
\author{Frank Sottile}
\address{\hskip-\parindent
        Mathematical Sciences Research Institute\\
        1000 Centennial Drive\\
        Berkeley, CA  \ 94720\\
        USA}

\curraddr{\hskip-\parindent
        Frank Sottile\\
        Department of Mathematics\\
        University of Wisconsin\\
        Van Vleck Hall\\
        480 Lincoln Drive\\
        Madison, Wisconsin 53706-1388\\
        USA}

\email{sottile@math.wisc.edu}
\urladdr{http://www.math.wisc.edu/\~{}sottile}
\date{1 April 1999}
\thanks{Research at MSRI supported by NSF grant DMS-9701755}

\begin{abstract}
We show that  the Schubert calculus of enumerative geometry is real, for
special Schubert conditions.
That is, for any such enumerative problem, there exist real conditions for
which all the {\it a priori} complex solutions are real.
\end{abstract}

\maketitle

\begin{center}
ERA of the AMS {\bf 5} (1999), pp.~35--39.
\end{center}
\bigskip

Fulton asked how many solutions to a problem of enumerative
geometry can be real, when that problem is one of counting geometric figures
of some kind having specified position with respect to some general fixed
figures~\cite{Fu_84}.
For the problem of plane conics tangent to five general conics, the
(surprising) answer is that all 3264 may be real~\cite{RTV}. 
Recently, Dietmaier has shown that all 40 positions of the Stewart
platform in robotics may be real~\cite{Dietmaier}.
Similarly, given any problem of enumerating lines in projective space
incident on some general fixed linear subspaces, there are real fixed
subspaces such that each of the (finitely many) incident lines are
real~\cite{Sottile97a}. 
Other examples are shown in~\cite{Sottile97c,Sottile97b},
and the case of 462 4-planes meeting 12 general 3-planes in ${\mathbb R}^7$
is due to an heroic symbolic computation~\cite{FRZ}.

For any problem of enumerating $p$-planes having excess  intersection with a
collection of fixed planes, we show there is a choice of fixed 
planes osculating a rational normal curve at real points so that each of the
resulting $p$-planes is real. 
This has implications for the problem of placing real poles in linear systems
theory~\cite{Byrnes} and is a special case of a far-reaching conjecture of
Shapiro and Shapiro~\cite{Sottile_shapiro}.

\section*{Special Schubert conditions}

For background on the Grassmannian, Schubert cycles, and the Schubert
calculus, see any of~\cite{Hodge_Pedoe,Griffiths_Harris,Fulton_tableaux}.
Let  $m,p\geq 1$ be integers.
Let $\gamma$ be a rational normal curve in ${\mathbb R}^{m+p}$.
For $k>0$ and $s\in\gamma$, let $K_k(s)$ be the $k$-plane
osculating $\gamma$ at $s$.
For every integer $a>0$, let $\tau_a(s)$ be the special Schubert cycle
consisting of $p$-planes $H$ in ${\mathbb C}^{m+p}$ 
which meet $K_{m+1-a}(s)$ nontrivially and let
$\tau^a(s)$ be the special Schubert cycle consisting of $p$-planes $H$
in ${\mathcal C}^{m+p}$ meeting $K_{m-1+a}(s)$ improperly: 
$\dim H\cap K_{m-1+a}(s)> a-1$.
These cycles $\tau_a(s)$ and $\tau^a(s)$ each have
codimension $a$ and $\tau^1=\tau_1$. 
Recall that the Grassmannian of $p$-planes in ${\mathbb C}^{m+p}$ has
dimension $mp$.
For any Schubert condition $w$, let $\sigma_w(s)$ be the Schubert cycle of
type $w$ given by the flag osculating $\gamma$ at $s$
and set $|w|$ to be the codimension of $\sigma_w(s)$.

\begin{thm}
Let $a_1,\ldots,a_n$ be positive integers with $a_1+\cdots+a_n=mp$.
For each $i=1,\ldots,n$ let $\sigma_i(s)$ be either $\tau^{a_i}(s)$ or
$\tau_{a_i}(s)$.
Then there exist real points $0,\infty,s_1,\ldots,s_n\in{\gamma}$ such that
for any Schubert conditions $w,v$, and integer $k$ with 
$|w|+|v|+a_k+\cdots+a_n=mp$, the intersection
\begin{equation}\label{eq:special}
  \sigma_w(0)\cap \sigma_v(\infty)\cap
  \sigma_k(s_k)\cap\cdots\cap\sigma_n(s_n)
\end{equation}
is transverse with all points of intersection real.
\end{thm}

Our proof is inspired by the Pieri homotopy algorithm of~\cite{HSS}.

We prove this in the case that each
$\sigma_i(s)=\tau_{a_i}(s)$.
This is no loss of generality, as the cycles $\tau_a(s)$ and $\tau^a(s)$
share the properties we need.
We use the following two results of Eisenbud and Harris~\cite{EH83},
who studied such intersections in their theory of limit linear systems.
For a Schubert class $w$, 
$w*a$ be the index of summation in the Pieri formula in the
cohomology of the Grassmannian~\cite{Fulton_tableaux},
$$
\sigma_w\cdot \tau_a\ =\  \sum_{v\in w*a} \sigma_v.
$$

\begin{prop}\label{prop:only}
\mbox{ }
\begin{enumerate}
\item (Theorem 2.3 of~\cite{EH83}) Let $s_1,\ldots,s_n$ be distinct points on
$\gamma$  and $w_1,\ldots,w_n$ be Schubert
conditions.
Then the intersection of Schubert cycles
$$
\sigma_{w_1}(s_1)\cap\sigma_{w_2}(s_2)\cap\cdots\cap\sigma_{w_n}(s_n)
$$
is proper in that it has dimension $mp-|w_1|-\cdots-|w_n|$.

\item (Theorem 8.1 of~\cite{EH83}) For any Schubert condition $w$,
integer $a>0$, and $0\in\gamma$,
we have 
$$
\lim_{t\rightarrow 0}\left( \sigma_w(0)\cap\tau_a(t)\rule{0pt}{11pt}\right)\
=\  \bigcup_{v\in w*a} \sigma_v(0),
$$
the limit taken along the rational normal curve, and as schemes.
\end{enumerate}
\end{prop}

\noindent{\bf Proof of Theorem~1.}
We argue by downward induction on $k$.
The initial case of $k=n$ holds as Pieri's formula implies the intersection
is a single, necessarily real, point.
Suppose it holds for $k$, and let $w,v$ satisfy 
$|w|+|v|+a_{k-1}+\cdots+a_n=mp$.

\noindent{\bf Claim}:   The cycle
$$
\sum_{u\in w*a_{k-1}}
\sigma_u(0)\cap\sigma_v(\infty)\cap
\tau_{a_k}(s_k)\cap\cdots\cap\tau_{a_n}(s_n)
$$
is free of multiplicities.

If not, then two summands, say $u$ and $u'$, have a point in common
and so
\begin{equation}\label{eq:bad_intersect}
  \sigma_u(0)\cap\sigma_{u'}(0)\cap\sigma_v(\infty)\cap
  \tau_{a_k}(s_k)\cap\cdots\cap\tau_{a_n}(s_n)
\end{equation}
is nonempty.
However, $\sigma_u(0)\cap\sigma_{u'}(0)$
is a Schubert cycle of smaller dimension.
Thus the intersection~(\ref{eq:bad_intersect}) 
must be empty, by Proposition~\ref{prop:only}~(1).
\medskip

From the claim and Proposition~\ref{prop:only}~(2), there is an
$\epsilon_{w,v}>0$ 
such that if $0<t\leq\epsilon_{w,v}$, then 
$$
\sigma_w(0)\cap \sigma_v(\infty)\cap
\tau_{a_{k-1}}(t)\cap\tau_{a_k}(s_k)\cap\cdots\cap\tau_{a_n}(s_n)
$$
is transverse with all points of intersection real.
Set 
$$
s_{k+1}\ :=\ min\{\epsilon_{w,v}\mid 
mp=|w|+|v|+a_{k-1}+\cdots+a_n\}.\qed
$$

\noindent{\bf Remark. }
Eisenbud and Harris~\cite{EH83} prove Proposition~\ref{prop:only}~(2) for
any nondegenerate arc $\gamma$.
Theorem~1 may be similarly strengthened.\medskip

\section*{Consequences}
Since small real perturbations of the points $s_1,\ldots,s_n$ cannot create
or destroy real points in a transverse intersection, Theorem~1 has the
following consequence.

\begin{cor}\label{cor:open}
There is an open subset (in the classical topology on $\gamma^n$) 
consisting of $n$-tuples $(s_1,\ldots,s_n)\in\gamma^n$ such
that~(\ref{eq:special}) is transverse with all points of intersection real.
\end{cor}

Theorem 1 proves part of a conjecture of Shapiro and Shapiro.

\begin{conj}[Shapiro and Shapiro]\label{conj:shapiro}
Let $m,p\geq 1$ and $w_1,\ldots,w_n$ be Schubert conditions on $p$-planes in
$(m+p)$-space. 
If $s_1,\ldots,s_n\in\gamma$ are real points so that 
\begin{equation}\label{eq:intersection}
\sigma_{w_1}(s_1)\cap\sigma_{w_2}(s_2)\cap\cdots\cap\sigma_{w_n}(s_n)
\end{equation}
is zero-dimensional, then all points of intersection are real.
\end{conj}

When the points $s_1,\ldots,s_n$ are are distinct and 
$|w_1|+\cdots+|w_n|=mp$,~(\ref{eq:intersection}) is
zero-dimensional~\cite{EH83}.
Besides Theorem~1, there is substantial and compelling evidence for the
validity of
Conjecture~\ref{conj:shapiro}~\cite{FRZ,Verschelde,RS98,Sottile_shapiro}. 
In every instance (choice of $w_1,\ldots,w_n$ and distinct real points 
$s_1,\ldots,s_n\in\gamma$) checked, all points in~(\ref{eq:intersection})
are real. 
This includes some when the $w_i$ are not special Schubert
conditions~\cite{Sottile_shapiro} and some spectacular
computations~\cite{FRZ,Verschelde}. 
For some sets of Schubert conditions, we have proven that for every choice of
distinct real points $s_1,\ldots,s_n\in\gamma$, all points
in~(\ref{eq:intersection}) are real. 
In every known case, the intersection scheme is reduced when the parameters
$s_1,\ldots,s_n$ are distinct and real. 
We conjecture this is always the case.

\begin{conj}\label{conj:reduced}
Let $m,p\geq 1$ and $w_1,\ldots,w_n$ be Schubert conditions with 
$mp=|w_1|+\cdots+|w_n|$. 
If $s_1,\ldots,s_n\in\gamma$ are distinct and real, then 
the intersection scheme~(\ref{eq:intersection}) is reduced.
\end{conj}

The number of real points in the scheme~(\ref{eq:intersection})
is locally constant on the set 
of parameters for which it is reduced.
Thus for special Schubert conditions, Conjecture~\ref{conj:reduced}
would imply Conjecture~\ref{conj:shapiro}, by Theorem 1.
In fact, more is true.

\begin{thm}
Conjecture~\ref{conj:reduced} implies Conjecture~\ref{conj:shapiro}.
\end{thm}

\noindent{\bf Proof. }
By Remark~3.4 to Theorem~3.3 of~\cite{Sottile_shapiro}, if
Conjecture~\ref{conj:shapiro} holds when $a_1=\cdots=a_{mp}=1$, then it
holds for any collection of Schubert conditions, if the
cycles~(\ref{eq:intersection}) are reduced for any choice of 
$w_1,\ldots,w_n$ with $mp=|w_1|+\cdots+|w_n|$ and general
$s_1,\ldots,s_n\in\gamma$. 
Both of these conditions are supplied by Conjecture~\ref{conj:reduced}.
\qed
\medskip

When $n=mp$ so that $a_1=\cdots=a_{mp}=1$,~(\ref{eq:special}) is a special
case of the pole placement problem in systems theory~\cite{Byrnes}. 
A physical system (e.g.~a mechanical linkage), called a {\it plant} with 
$m$ inputs and $p$ measured outputs 
whose evolution is governed by a system of linear differential equations
is modeled by a
$m\times(m+p)$-matrix of univariate polynomials $[N(s): D(s)]$.
The largest degree of a maximal minor of this matrix is the MacMillan
degree, $n$, of the evolution equation.
For $s_1,\ldots,s_n\in{\mathbb C}$, any $p$-plane $H$ satisfying
$$
H\cap \mbox{\rm row span}[N(s_i): D(s_i)]\ \neq\ \{0\}\qquad
i=1,\ldots,n
$$
gives a constant linear feedback law for which the resulting closed system
has natural frequencies (poles) $s_1,\ldots,s_n$.

In this way, the pencil $K(s)$ of $m$-planes osculating $\gamma$ 
represents a particular $m$-input $p$-output plant of McMillan degree
$mp$, and the $p$-planes in
$$
\tau_1(s_1)\cap\tau_1(s_2)\cap\cdots\cap\tau_1(s_{mp})
$$
represent
linear feedback laws for which the resulting closed system has
natural frequencies $s_1,\ldots,s_{mp}$.
Since translation along $\gamma$ fixes the system and acts on the feedback
laws by a real linear transformation, we may assume that
$s_1,\ldots,s_{mp}<0$ and so we obtain the existence of a stable 
system with only real feedback laws.

\begin{cor}
The plant represented by $K(s)$ can be stabilized by real feedback laws in
such a way that all feedback laws are real.
\end{cor}

\end{document}